\documentclass[titlepage,draft,12pt]{article} 
\usepackage{amsfonts,latexsym,pstricks,pst-plot} 
\usepackage[a4paper]{geometry}

\newcommand{\ocont}{$\omega$-continuous}
\newcommand{\G}{\mathcal{G}}
\newcommand{\GM}{\mathcal{GM}}
\newcommand{\trebar}[1]{\|#1\|}
\newcommand{\lk}{\lambda_{k}}
\newcommand{\ep}{\varepsilon}
\newcommand{\re}{{\mathbb{R}}}
\newcommand{\n}{{\mathbb{N}}}

\newcommand{\vs}{\vspace{.2cm}}
\newcommand{\prf}{{\sc Proof.}\ }
\newcommand{\qed}{{\penalty 10000\mbox{$\quad\Box$}}\bigskip}

\newcommand{\m}[1]{m(\au{#1}^{2})}
\newcommand{\au}[1]{|A^{1/2}#1|}

\newcommand{\cep}{c_{\ep}}

\newtheorem{thm}{Theorem}[section]
\newtheorem{thmbibl}{Theorem}
\newtheorem{rmk}[thm]{Remark}
\newtheorem{prop}[thm]{Proposition}

\newtheorem{lemma}[thm]{Lemma}

\title{Spectral gap global solutions for degenerate Kirchhoff equations}
\author{Marina Ghisi\vs\\ {\normalsize
Universit\`a degli Studi di Pisa} \\{\normalsize Dipartimento di
Matematica ``Leonida Tonelli''}\\
{\normalsize 
PISA (Italy)}\\  
{\normalsize e-mail: \texttt{ghisi@dm.unipi.it}}\and 
Massimo Gobbino\vs\\ {\normalsize Universit\`a degli Studi di Pisa} 
\\{\normalsize Dipartimento di Matematica Applicata ``Ulisse Dini''}\\ 
{\normalsize 
 PISA (Italy)}\\  
{\normalsize e-mail: \texttt{m.gobbino@dma.unipi.it}}}
\date{}

\begin{document}
\maketitle
\begin{abstract}
	We consider the second order Cauchy problem 
	$$u''+\m{u}Au=0, 
	\hspace{2em}
	u(0)=u_{0},\quad
	u'(0)=u_{1},$$
	where $m:[0,+\infty)\to[0,+\infty)$ is a continuous function, and
	$A$ is a self-adjoint nonnegative operator with dense domain on a
	Hilbert space.
	
	It is well known that this problem admits local-in-time solutions
	provided that $u_{0}$ and $u_{1}$ are regular enough, depending on
	the continuity modulus of $m$, and on the strict/weak
	hyperbolicity of the equation. 
	
	We prove that for such initial data $(u_{0},u_{1})$ there exist
	two pairs of initial data $(\overline{u}_{0},\overline{u}_{1})$,
	$(\widehat{u}_{0},\widehat{u}_{1})$ for which the solution is
	global, and such that $u_{0}=\overline{u}_{0}+\widehat{u}_{0}$,
	$u_{1}=\overline{u}_{1}+\widehat{u}_{1}$.
	
	This is a byproduct of a global existence result for initial data 
	with a suitable spectral gap, which extends previous results
	obtained in the strictly hyperbolic case with a smooth
	nonlinearity $m$.
	
\vspace{1cm}

\noindent{\bf Mathematics Subject Classification 2000 (MSC2000):}
35L70, 35L80, 35L90.

\vspace{1cm} 

\end{abstract}
 
\section{Introduction}

Let $H$ be a real Hilbert space.  For every $x$ and $y$ in $H$, let
$|x|$ denote the norm of $x$, and let $\langle x,y\rangle$ denote the
scalar product of $x$ and $y$.  Let $A$ be an unbounded linear
operator on $H$ with dense domain $D(A)$.  We always assume that $A$
is self-adjoint and nonnegative, so that for every $\alpha\geq 0$ the
power $A^{\alpha}$ is defined in a suitable domain $D(A^{\alpha})$.

Given a continuous function $m:[0,+\infty)\to[0,+\infty)$ we consider 
the Cauchy problem
\begin{equation}
	u''(t)+\m{u(t)}Au(t)=0, 
	\hspace{2em}\forall t\in[0,T),
	\label{pbm:h-eq}
\end{equation}
\begin{equation}
	u(0)=u_0,\hspace{3em}u'(0)=u_1.
	\label{pbm:h-data}
\end{equation}
 
It is well known that (\ref{pbm:h-eq}), (\ref{pbm:h-data}) is the
abstract setting of the Cauchy-boundary value problem for the
quasilinear hyperbolic integro-differential partial differential
equation
\begin{equation}
	u_{tt}(t,x)-
	m{\left(\int_{\Omega}\left|\nabla u(t,x)\right|^2\,dx\right)}
	\Delta u(t,x)=0
	\hspace{2em}
	\forall(x,t)\in\Omega\times[0,T),
	\label{eq:k}
\end{equation}
where $\Omega\subseteq\re^{n}$ is an open set, and $\nabla u$ and
$\Delta u$ denote the gradient and the Laplacian of $u$ with respect
to the space variables.

Equation (\ref{pbm:h-eq}) is called strictly hyperbolic if
\begin{equation}
	m(\sigma)\geq\nu>0
	\quad\quad
	\forall\sigma\geq 0.
	\label{hp:s-h}
\end{equation}

Equation (\ref{pbm:h-eq}) is called weakly (or degenerate) hyperbolic
if 
$$m(\sigma)\geq 0 \quad\quad
\forall\sigma\geq 0.$$

Existence of local/global solutions to (\ref{pbm:h-eq}),
(\ref{pbm:h-data}) has long been investigated in the last century.
The theory is well established in the case of local solutions, which
are known to exist in the following situations.
\begin{enumerate}
	\item[(L1)] When equation (\ref{pbm:h-eq}) is strictly hyperbolic,
	$m$ is Lipschitz continuous, and initial data $(u_{0},u_{1})\in
	D(A^{3/4})\times D(A^{1/4})$ (see \cite{ap} and the references
	quoted therein).

	\item[(L2)] When equation (\ref{pbm:h-eq}) is weakly hyperbolic,
	$m$ is continuous, and initial data are analytic.  In this case
	solutions are actually global (see \cite{as}, \cite{das-an-1},
	\cite{das-an-2}).

	\item[(L3)] More generally, when initial data belong to suitable
	intermediate spaces, depending on the continuity modulus of $m$,
	and on the strict/weak hyperbolicity of (\ref{pbm:h-eq}) (see
	\cite{hirosawa-main} and \cite{gg:k-derloss}).  This is a sort of
	interpolation between (L1) and (L2).  We refer to
	section~\ref{sec:prelim} for precise definitions of the functional
	spaces in the abstract framework and a formal local existence
	statement (Theorem~\ref{thm:hirosawa}).
\end{enumerate}

Existence of global solutions is a much more difficult problem, and it
is still widely open.  A positive answer has been given in the case
(L2), and in some special situations: quasi-analytic initial data (see
\cite{nishihara}), or Sobolev-type data but special nonlinearities $m$
(see \cite{poho-m}), or dispersive operators and small data (see
\cite{gh}, \cite{das}).  But for (L2) all these results assume the
strict hyperbolicity and the Lipschitz continuity of $m$.

Recently \textsc{R.\ Manfrin}~\cite{manfrin1,manfrin2} (see
also~\cite{hirosawa2}) considered once again the strictly hyperbolic
case with a smooth nonlinearity.  He proved global existence in a
special class of nonanalytic initial data.  Manfrin's spaces are not
vector spaces and do not contain any Gevrey space $\G_{s}$ with $s>1$.
However they have the following astonishing property: 
\begin{enumerate}
	\item[(M)] every pair of initial conditions $(u_{0},u_{1})\in
	D(A)\times D(A^{1/2})$ is the sum of two pairs of initial
	conditions in Manfrin's spaces, i.e., the sum of two initial
	conditions for which the solution is global!
\end{enumerate}

This theory requires the strict hyperbolicity and some smoothness of
$m$, which is assumed to be of class $C^{2}$ both in \cite{manfrin1}
and \cite{manfrin2}. 

In this paper we extend Manfrin's theory to the general situation of
(L3).  We consider indeed both the strictly hyperbolic and the weakly
hyperbolic case, and a nonlinearity $m$ with a given continuity
modulus.  In Theorem~\ref{thm:main} we prove global existence for
initial data in a suitable subset of the spaces involved in (L3).  In
analogy with Manfrin's spaces, the definition (\ref{defn:m-space}) of
our subset is made in terms of the spectral resolution of initial
data.  Of course our subset is not a vector space and it doesn't even
contain all analytic functions.  Nevertheless in
Proposition~\ref{prop:sum} we show that this subset satisfies property
(M) in the spaces involved in (L3).

From the point of view of property (M) our result extends Manfrin's
one also in the framework (L1).  In this case we obtain indeed
property (M) for initial data in $D(A^{3/4})\times D(A^{1/4})$ and a
locally Lipschitz continuous nonlinearity $m$, instead of initial data
in $D(A)\times D(A^{1/2})$ and $m\in C^{2}$.

This paper is organized as follows.  In section~\ref{sec:prelim} we
recall the definition of continuity modulus and Gevrey-type functional
spaces, and we state the local existence result for the case (L3).  In
section~\ref{sec:statements} we introduce our spaces and we state our
main results.  In section~\ref{sec:proofs} we prove these results.

\setcounter{equation}{0}
\section{Preliminaries}\label{sec:prelim}

For the sake of simplicity we assume that $H$ admits a countable
complete orthonormal system $\{e_{k}\}_{k\geq 1}$ made by eigenvectors
of $A$.  We denote the corresponding eigenvalues by $\lambda_{k}^{2}$
(with $\lk\geq 0$), so that $Ae_{k}=\lambda_{k}^{2}e_{k}$ for every
$k\geq 1$.

Under this assumption we can work with Fourier series.  However, any
definition or statement of this section can be easily extended to the
general setting just by using the spectral decomposition instead of
Fourier series.  The interested reader is referred to \cite{ap} for
further details.

By means of the orthonormal system every $u\in H$ can be written in a
unique way in the form $u=\sum_{k=1}^{\infty}u_{k}e_{k}$, where
$u_{k}=\langle u,e_{k}\rangle$ are the Fourier components of $u$.
With these notations for every $\alpha\geq 0$ we have that 
$$D(A^{\alpha}):=\left\{u\in H:\sum_{k=1}^{\infty}
\lambda_{k}^{4\alpha}u_{k}^{2}<+\infty\right\}.$$

Let now $\varphi:[0,+\infty)\to(0,+\infty)$ be any function.  Then for
every $\alpha\geq 0$ and $r>0$ one can set
\begin{equation}
	\trebar{u}_{\varphi,r,\alpha}^{2}:=\sum_{k=1}^{\infty}\lambda_{k}^{4\alpha}
	u_{k}^{2} \exp\left(\strut r\varphi(\lambda_{k})\right),
	\label{defn:trebar}
\end{equation}
and then define the spaces
$$\G_{\varphi,r,\alpha}(A):=
\left\{u\in H:\trebar{u}_{\varphi,r,\alpha}^{2}<+\infty\right\}.$$

These spaces are a generalization of the usual spaces of Sobolev,
Gevrey or analytic functions.  They are Hilbert spaces with norm
$(|u|^{2}+\trebar{u}_{\varphi,r,\alpha}^{2})^{1/2}$. We also set
$$\G_{\varphi,\infty,\alpha}(A):=\bigcap_{r>0}\G_{\varphi,r,\alpha}(A).$$

A \emph{continuity modulus} is a continuous increasing function
$\omega:[0,+\infty)\to[0,+\infty)$ such that $\omega(0)=0$, and
$\omega(a+b)\leq\omega(a)+\omega(b)$ for every $a\geq 0$ and $b\geq
0$.

The function $m$ is said to be $\omega$-continuous if there exists a
constant $L\in\re$ such that
\begin{equation}
	|m(a)-m(b)|\leq
	L\,\omega(|a-b|)
	\hspace{3em}
	\forall a\geq 0,\ \forall b\geq 0.
	\label{hp:m-ocont}
\end{equation}

The following result sums up the state of the art concerning existence
of local solutions.  We refer to Theorem~2.1 and Theorem~2.2
in~\cite{hirosawa-main} for the existence part, to \cite{gg:k-derloss}
for some counterexamples, and to \cite{gg:k-uniq} for uniqueness
issues.

\begin{thmbibl}\label{thm:hirosawa}
	Let $\omega$ be a continuity modulus, let
	$m:[0,+\infty)\to[0,+\infty)$ be a (locally) \ocont\ function, and
	let $\varphi:[0,+\infty)\to(0,+\infty)$.
	
	Let us assume that there exists a constant $\Lambda$ such that
	\begin{equation}
		\sigma
		\omega\left(\frac{1}{\sigma}\right)\leq\Lambda\varphi(\sigma)
		\quad\quad
		\forall\sigma> 0
		\label{hp:phi-ndg}
	\end{equation}
	in the strictly hyperbolic case, and 
	\begin{equation}
		\sigma\leq\Lambda\varphi\left(\frac{\sigma}{
		\sqrt{\omega(1/\sigma)}}\right)
		\quad\quad
		\forall\sigma> 0
		\label{hp:phi-dg}
	\end{equation}
	in the weakly hyperbolic case.
	
	Let 
	\begin{equation}
		(u_{0},u_{1})\in
			\G_{\varphi,r_{0},3/4}(A)\times\G_{\varphi,r_{0},1/4}(A)
		\label{hp:hiro-data}
	\end{equation}
	for some $r_{0}>0$.
	
	Then there exists $T>0$, and a nonincreasing function
	$r:[0,T]\to(0,r_{0}]$ such that problem (\ref{pbm:h-eq}),
	(\ref{pbm:h-data}) admits at least one local solution 
	\begin{equation}
		u\in C^{1}\left([0,T];\G_{\varphi,r(t),1/4}(A)\right)\cap
		C^{0}\left([0,T];\G_{\varphi,r(t),3/4}(A)\right).  
	\label{th:reg-sol}
\end{equation}
\end{thmbibl}

\setcounter{equation}{0}
\section{Main result}\label{sec:statements}

Let $\mathcal{L}$ denote the set of all sequences $\{\rho_{n}\}$ of
positive real numbers such that $\rho_{n}\to +\infty$ as $n\to
+\infty$.  Given $\varphi:[0,+\infty)\to(0,+\infty)$,
$\{\rho_{n}\}\in\mathcal{L}$, $\alpha\geq 0$, and $\beta\geq 0$ we
set
\begin{equation}
	\GM_{\varphi,\{\rho_{n}\},\alpha}^{(\beta)}(A):=\left\{ u\in
	H:\sum_{\lk>\rho_{n}}\lambda_{k}^{4\alpha}u_{k}^{2}
	\exp\left(\rho_{n}^{\beta}\varphi(\lk)\right)\leq\rho_{n}
	\quad\forall n\in\n\right\},
	\label{defn:m-space}
\end{equation}
and then
$$\GM_{\varphi,\alpha}^{(\beta)}(A):=
\bigcup_{\{\rho_{n}\}\in\mathcal{L}}
\GM_{\varphi,\{\rho_{n}\},\alpha}^{(\beta)}(A).$$

These spaces are a generalization of Manfrin's spaces.

The following global existence result is the main result of this
paper.

\begin{thm}\label{thm:main}
	Let $\omega$ be a continuity modulus, let
	$m:[0,+\infty)\to[0,+\infty)$ be a function satisfying
	(\ref{hp:m-ocont}), let $\varphi:[0,+\infty)\to(0,+\infty)$, and
	let $\{\rho_{n}\}\in\mathcal{L}$.
	
	Let us assume that
	\begin{itemize}
		\item  in the strictly hyperbolic case (\ref{hp:phi-ndg})
		holds true for a suitable $\Lambda$, and
		\begin{equation}
			(u_{0},u_{1})\in
			\GM_{\varphi,\{\rho_{n}\},3/4}^{(2)}(A)\times
			\GM_{\varphi,\{\rho_{n}\},1/4}^{(2)}(A),
			\label{hp:data-ndg}
		\end{equation}  
	
		\item  in the weakly hyperbolic case (\ref{hp:phi-dg})
		holds true for a suitable $\Lambda$, and
		\begin{equation}
			(u_{0},u_{1})\in
			\GM_{\varphi,\{\rho_{n}\},3/4}^{(3)}(A)\times
			\GM_{\varphi,\{\rho_{n}\},1/4}^{(3)}(A).
			\label{hp:data-dg}
		\end{equation}  
	\end{itemize}
	
	Then problem (\ref{pbm:h-eq}), (\ref{pbm:h-data})
	admits at least one global solution $u(t)$ with
	\begin{equation}
		u\in C^{1}\left([0,+\infty);\G_{\varphi,r,3/4}(A)\right)\cap
		C^{0}\left([0,+\infty);\G_{\varphi,r,1/4}(A)\right)
		\label{hp:reg-sol}
	\end{equation}
	for every $r>0$.
\end{thm}

We conclude by speculating on these spaces. First of all it is easy to
prove that
\begin{equation}
	\GM_{\varphi,\alpha}^{(\beta)}(A)\subseteq
	\G_{\varphi,\infty,\alpha}(A)
	\label{eq:inclusion}
\end{equation}
for every admissible values of the parameters.  On one hand this
inclusion is ``very strict''.  Roughly speaking indeed the
inequalities in definition (\ref{defn:m-space}) require that the
spectrum of $u$ ``has a big hole after each $\rho_{n}$''.  For this
heuristic reason we used ``spectral gap solutions'' to denote the
solutions produced by Theorem~\ref{thm:main}.

On the other hand inclusion (\ref{eq:inclusion}) is ``not so strict''
in the sense that
$$\GM_{\varphi,\alpha}^{(\beta)}(A)+\GM_{\varphi,\alpha}^{(\beta)}(A)=
\G_{\varphi,\infty,\alpha}(A)$$
for any admissible values of the parameters.  We state this property
more precisely in the case of pairs of initial data.

\begin{prop}\label{prop:sum}
	Let $\varphi:[0,+\infty)\to(0,+\infty)$, and let
	\begin{equation}
		(u_{0},u_{1})\in
		\G_{\varphi,\infty,3/4}(A)\times\G_{\varphi,\infty,1/4}(A).
		\label{hp:prop}
	\end{equation}
	
	Then for every $\beta\geq 0$ there exist $\{\overline{\rho}_{n}\}$ and 
	$\{\widehat{\rho}_{n}\}$ in $\mathcal{L}$, and
	\begin{equation}
		(\overline{u}_{0},\overline{u}_{1})\in
		\GM_{\varphi,\{\overline{\rho}_{n}\},3/4}^{(\beta)}(A)\times
		\GM_{\varphi,\{\overline{\rho}_{n}\},1/4}^{(\beta)}(A),
	\label{th:prop-1}
\end{equation}
	\begin{equation}
		(\widehat{u}_{0},\widehat{u}_{1})\in
		\GM_{\varphi,\{\widehat{\rho}_{n}\},3/4}^{(\beta)}(A)\times
		\GM_{\varphi,\{\widehat{\rho}_{n}\},1/4}^{(\beta)}(A),
		\label{th:prop-2}
	\end{equation}
	such that $u_{0}=\overline{u}_{0}+\widehat{u}_{0}$ and
	$u_{1}=\overline{u}_{1}+\widehat{u}_{1}$.
\end{prop}

\begin{rmk}
	\begin{em}		
		Combining Theorem~\ref{thm:main} and
		Proposition~\ref{prop:sum} we obtain the following statement:
		every pair of initial conditions satisfying
		(\ref{hp:hiro-data}) with $r_{0}=\infty$ is the sum of two
		pairs of initial conditions for which the solution is global.
		We have thus extended to the general case the astonishing
		aspect of Manfrin's result.
		
		The extra requirement that $r_{0}=\infty$ is hardly
		surprising.  It is indeed a necessary condition for existence
		of global solutions even in the theory of linear equations
		with nonsmooth time dependent coefficients.
	\end{em}
\end{rmk}

\begin{rmk}
	\begin{em}
		The $\omega$-continuity assumption on $m$ can be easily
		relaxed to local $\omega$-continuity in all the cases where
		there is a uniform-in-time estimate of $|A^{1/2}u(t)|$ in
		terms of the initial data.  We refer to the paragraph ``Energy
		conservation'' in section~\ref{sec:proof-prelim} for further
		details.
	\end{em}
\end{rmk}
		
\begin{rmk}\label{rmk:reg-sol}
	\begin{em}
		It is possible to extend the result of Theorem~\ref{thm:main}
		to larger spaces.  A careful inspection of the proof reveals
		that in the strictly hyperbolic case one can replace $\beta=2$
		with any $\beta>1$, in the weakly hyperbolic case one can
		replace $\beta=3$ with any $\beta>2$.  It should also be
		possible to enlarge these spaces in order to contain all
		analytic functions, for which a global solution was already
		known to exist.
		
		Our choice (\ref{defn:m-space}) is optimized in order to
		obtain both Theorem~\ref{thm:main} and
		Proposition~\ref{prop:sum} under the more general assumptions
		on $m$, and with a simple proof.
	\end{em}
\end{rmk}

\setcounter{equation}{0}
\section{Proofs}\label{sec:proofs}

\subsection{Preliminaries}\label{sec:proof-prelim}

\paragraph{Estimates for a continuity modulus}

The following estimates are crucial in the proof of our main result
(see also Lemma~3.1 in \cite{gg:k-uniq}).

\begin{lemma}
	Let $\omega:[0,+\infty)\to[0,+\infty)$ be a continuity modulus.
	
	Then
	\begin{eqnarray}
		 & \omega(\lambda x)\leq(1+\lambda)\omega(x)
		 \quad\quad\forall\lambda\geq 0,\ \forall x\geq 0; & 
		\label{th:omega-lambda}  \\
		\noalign{\vspace{1ex}}
		 & \displaystyle{\omega(x)\geq\omega(1)\frac{x}{x+1}}
		 \quad\quad \forall x\geq 0; & 
		\label{th:omega-est}  \\
		\noalign{\vspace{1ex}}
		 & \displaystyle{1+\frac{1}{\omega(x)}\leq\left(1+\frac{1}{\omega(1)}\right)
		 \left(1+\frac{1}{x}\right)}
		 \quad\quad\forall x>0.& 
		\label{th:omega-3}
	\end{eqnarray}
\end{lemma}

\prf
Inequality (\ref{th:omega-lambda}) can be easily proved by induction
on the integer part of $\lambda$ using the monotonicity and the
subadditivity of $\omega$.  Inequality (\ref{th:omega-est}) follows
from (\ref{th:omega-lambda}) applied with $\lambda=1/x$.  Inequality
(\ref{th:omega-3}) follows from (\ref{th:omega-est}).
\qed

\paragraph{Energy conservation}

Let $u$ be any solution of (\ref{pbm:h-eq}) defined in an interval
$[0,T)$. Let us set
$$M(\sigma):=\int_{0}^{\sigma}m(s)\,ds
\quad\quad
\forall\sigma\geq 0,$$
and let us consider the usual Hamiltonian
$$\mathcal{H}(t):=|u'(t)|^{2}+M(|A^{1/2}u(t)|^{2}).$$

It is well known that $\mathcal{H}(t)$ is constant. In particular
\begin{equation}
	|u'(t)|^{2}\leq\mathcal{H}(0)
	\quad\quad
	\forall t\in[0,T).
	\label{est:u'}
\end{equation}

In the strictly hyperbolic case we have also that
$M(\sigma)\geq\nu\sigma$, hence
\begin{equation}
	|A^{1/2}u(t)|^{2}\leq\frac{\mathcal{H}(0)}{\nu}
	\quad\quad
	\forall t\in[0,T).
	\label{est:au}
\end{equation}

This provides an estimate of $|A^{1/2}u(t)|$ in terms of the initial
conditions.  This type of estimate can be obtained also without the
strict hyperbolicity provided that the limit of $M(\sigma)$ as
$\sigma\to +\infty$ is $+\infty$ or at least larger than
$\mathcal{H}(0)$.

\paragraph{Convolutions}

In the next result we recall the properties of convolutions which are 
needed in the sequel (we omit the standard proof).

\begin{lemma}\label{lemma:conv}
	Let $\rho:\re\to [0,+\infty)$ be a function of class $C^{\infty}$,
	with support contained in $[-1,1]$, and integral equal to 1.
	
	Let $a>0$, and let $f:[0,a]\to\re$ be a continuous function.  Let
	us extend $f(x)$ to the whole real line by setting $f(x)=f(0)$ for
	every $x\leq 0$, and $f(x)=f(a)$ for every $x\geq a$.  
	
	For every $\ep>0$ let us set 
	$$f_{\ep}(x):=\int_{\re}^{}f(x+\ep
	s)\rho(s)\,ds \quad\quad
	\forall x\in\re.$$
	
	Then $f_{\ep}(x)$ has the following properties.
	\begin{enumerate}
		\renewcommand{\labelenumi}{(\arabic{enumi})}
		\item  $f_{\ep}\in C^{\infty}(\re)$.
	
		\item If $\mu_{1}\leq f(x)\leq\mu_{2}$ for every $x\in[0,a]$,
		then $\mu_{1}\leq f_{\ep}(x)\leq\mu_{2}$ for every
		$x\in\re$ and every $\ep>0$.
	
		\item  $|f_{\ep}(0)|\leq\max\{|f(x)|:0\leq x\leq\ep\}$ for
		every $\ep> 0$.
	
		\item Let $\omega$ be a continuity modulus. Let us assume that 
		\begin{equation}
			|f(x)-f(y)|\leq H\omega(|x-y|) \quad\quad \forall x\in[0,a],\
			\forall y\in[0,a],
			\label{hp:o-cont}
		\end{equation}
		for some $H\geq 0$.  Then there exists a constant
		$\gamma_{0}$ (independent on $\ep$, $H$, and on the function
		$f(t)$) such that
		$$|f_{\ep}(x)-f(x)|\leq\gamma_{0}H\omega(\ep)
			\quad\quad
			\forall x\in\re,\ \forall\ep>0,$$
		$$|f_{\ep}'(x)|\leq\gamma_{0}H\,\displaystyle{\frac{\omega(\ep)}{\ep}}
			\quad\quad
			\forall x\in\re,\ \forall\ep>0.$$

	\end{enumerate}
\end{lemma}

\paragraph{Maximal local solutions}

By (\ref{eq:inclusion}) assumptions (\ref{hp:data-ndg}) and
(\ref{hp:data-dg}) imply that
$(u_{0},u_{1})\in\G_{\varphi,\infty,3/4}(A)\times
\G_{\varphi,\infty,1/4}(A)$.  Therefore the existence of a local
solution to (\ref{pbm:h-eq}), (\ref{pbm:h-data}) follows from
Theorem~\ref{thm:hirosawa} both in the strictly hyperbolic and in the
weakly hyperbolic case.  Since initial data satisfy
(\ref{hp:hiro-data}) for every $r_{0}$, from the linear theory it
easily follows that the local solution satisfies (\ref{th:reg-sol})
for every $r(t)$.

By a standard argument any local solution can be continued to a
solution defined in a maximal interval $[0,T)$.  If $T=+\infty$ there
is nothing to prove.  In order to exclude that $T<+\infty$ we prove
that the time derivative of $|A^{1/2}u(t)|^{2}$ cannot blow-up in a
finite time.  The proof of this a priori estimate, which is the basic
tool in all global existence results, is different in the strictly
hyperbolic and in the weakly hyperbolic case.

\subsection{The strictly hyperbolic case}

Let us introduce some constants. From the strict hyperbolicity
(\ref{hp:s-h}) and estimate (\ref{est:au}) we have that
$$\nu\leq\m{u(t)}\leq
\max\left\{m(\sigma):0\leq\sigma\leq
\frac{\mathcal{H}(0)}{\nu}\right\}=:\mu
\quad\quad
\forall t\geq 0.$$

Let $L$, $\Lambda$, $\gamma_{0}$ be the constants appearing in
(\ref{hp:m-ocont}), (\ref{hp:phi-ndg}), and in Lemma~\ref{lemma:conv},
and let
$$\gamma_{1}:=\max\{1,\mu\}\cdot\max\left\{1,\nu^{-1}\right\},$$
$$H_{1}:=\max\left\{\left|\langle A^{3/4}u_{0},A^{1/4}u_{1}
\rangle\right|+1,\left(1+\nu^{-1}\right)
\mathcal{H}(0)+2\gamma_{1}+1\right\},$$
$$\gamma_{2}:=\gamma_{0}L\Lambda(2H_{1}+1)\left(
\frac{1}{\nu}+\frac{1}{\sqrt{\nu}}\right).$$

Since $\rho_{n}\to +\infty$ we can choose $n\in\n$ such that
\begin{equation}
	\rho_{n}\geq\max\{\gamma_{2}T,1\}.
	\label{defn:rhon}
\end{equation}

Let us set 
$$S:=\sup\left\{\tau\leq T:\left|\langle A^{3/4}u(t),A^{1/4}u'(t)
\rangle\right|\leq H_{1}\rho_{n}\;\;\forall t\in[0,\tau]\right\}.$$

We remark that $S>0$ because
$\left|\langle A^{3/4}u_{0},A^{1/4}u_{1}
\rangle\right|<H_{1}\leq H_{1}\rho_{n}$.

Now we distinguish the case $S=T$ and $S<T$.

\subparagraph{\textmd{\emph{Case}} $S=T$}

The argument is quite standard. In the interval $[0,T)$ the function
$u(t)$ is the solution of the linear problem
\begin{equation}
	v''(t)+c(t)Av(t)=0
	\label{pbm:lin-eq}
\end{equation}
\begin{equation}
	v(0)=u_{0},
	\quad\quad
	v'(0)=u_{1},
	\label{pbm:lin-data}
\end{equation}
where
\begin{equation}
	c(t):=\m{u(t)}.
	\label{defn:c}
\end{equation}

Since $S=T$ in this case we have that
\begin{equation}
	\left|\frac{\mathrm{d}}{\mathrm{d}t}|A^{1/2}u(t)|^{2}\right|=
	2\left|\langle A^{3/4}u(t),A^{1/4}u'(t)
	\rangle\right|\leq 2H_{1}\rho_{n}
	\label{eq:c'}
\end{equation}
for every $t\in[0,T)$. It follows that $|A^{1/2}u(t)|^{2}$ is Lipschitz
continuous in $[0,T)$, hence $c(t)$ can be extended to an
$\omega$-continuous function defined in the closed interval $[0,T]$.
By the linear theory (see \cite{dgcs} and \cite{hirosawa-main})
problem (\ref{pbm:lin-eq}), (\ref{pbm:lin-data}) has a solution
$$v\in C^{0}\left([0,T];\G_{\varphi,r,3/4}(A)\right)\cap
C^{1}\left([0,T];\G_{\varphi,r,1/4}(A)\right)$$
for every $r>0$. Since the solution of the linear problem is unique,
this implies that there exist
$$\widehat{u}_{0}:=\lim_{t\to T^{-}}u(t)\in\G_{\varphi,\infty,3/4}(A),
\quad\quad
\widehat{u}_{1}:=\lim_{t\to T^{-}}u'(t)\in\G_{\varphi,\infty,1/4}(A).$$

Applying Theorem~\ref{thm:hirosawa} with initial data
$(\widehat{u}_{0},\widehat{u}_{1})$ one can therefore continue $u(t)$
on an interval $[0,T_{1})$ with $T_{1}>T$, which contradicts the
maximality of $T$.

\subparagraph{\textmd{\emph{Case}} $S<T$}

By the maximality of $S$ we have that necessarily
\begin{equation}
	\left|\langle A^{3/4}u(S),A^{1/4}u'(S) \rangle\right|=
	H_{1}\rho_{n}.
	\label{eq:S-nec}
\end{equation}

Let us consider the function $c(t)$ defined according to
(\ref{defn:c}).  In this case (\ref{eq:c'}) holds true for every $t\in[0,S]$,
hence by (\ref{hp:m-ocont}) and (\ref{th:omega-lambda}) we have that
\begin{eqnarray*}
	\left|c(t)-c(s)\right| & = & \left|\m{u(t)}-\m{u(s)}\right|  \\
	 & \leq & L\,\omega\left(\left|
	 |A^{1/2}u(t)|^{2}-|A^{1/2}u(s)|^{2}\right|\right)  \\
	 & \leq & L\,\omega(2H_{1}\rho_{n}|t-s|)  \\
	 & \leq & L(2H_{1}\rho_{n}+1)\,\omega(|t-s|)  \\
	 & \leq & L(2H_{1}+1)\rho_{n}\,\omega(|t-s|) 
\end{eqnarray*}
for every $t$ and $s$ in $[0,S]$.  Let us extend $c(t)$ outside the
interval $[0,S]$ as in Lemma~\ref{lemma:conv}, and let us set
\begin{equation}
	\cep(t):=\int_{\re}^{}c(t+\ep s)\rho(s)\,ds \quad\quad \forall
	t\in\re.
	\label{defn:cep}
\end{equation}

Since estimate (\ref{hp:o-cont}) holds true with
$H:=L(2H_{1}+1)\rho_{n}$, from statements (2) and (4) of
Lemma~\ref{lemma:conv} we deduce that
\begin{equation}
	\nu\leq\cep(t)\leq\mu
	\quad\quad
	\forall t\in\re,\ \forall\ep>0,  
	\label{est:cep-nu-mu}
\end{equation}
\begin{equation}
	 |\cep(t)-c(t)|\leq\gamma_{0}L(2H_{1}+1)\rho_{n}\omega(\ep)
	\quad\quad
	\forall t\in\re,\ \forall\ep>0,  
	\label{est:cep-c}
\end{equation}
\begin{equation}
	|\cep'(t)|\leq\gamma_{0}L(2H_{1}+1)\rho_{n}\,
	\displaystyle{\frac{\omega(\ep)}{\ep}}
	\quad\quad
	\forall t\in\re,\ \forall\ep>0. 
	\label{est:cep'}
\end{equation}

Let us consider the Fourier components $u_{k}(t)$ of $u(t)$, and let
us set 
\begin{equation}
	E_{k,\ep}(t):=|u_{k}'(t)|^{2}+\lk^{2}\cep(t)|u_{k}(t)|^{2}.
	\label{defn:ekep}
\end{equation}

An easy computation shows that
\begin{eqnarray*}
	E_{k,\ep}'(t) & = & \cep'(t)\lk^{2}|u_{k}(t)|^{2}+
	2\lk^{2}(\cep(t)-c(t))u_{k}(t)u_{k}'(t)\\
	 & \leq & \frac{|\cep'(t)|}{\cep(t)}\cep(t)\lk^{2}|u_{k}(t)|^{2}+
	 \lk\frac{|\cep(t)-c(t)|}{\sqrt{\cep(t)}}2|u_{k}'(t)|\cdot
	 \lk\sqrt{\cep(t)}|u_{k}(t)|\\
	 & \leq & \frac{|\cep'(t)|}{\cep(t)}E_{k,\ep}(t)+
	 \lk\frac{|\cep(t)-c(t)|}{\sqrt{\cep(t)}}E_{k,\ep}(t),
\end{eqnarray*}
hence by (\ref{est:cep-nu-mu}), (\ref{est:cep-c}), and
(\ref{est:cep'}) we obtain that
\begin{equation}
	E_{k,\ep}'(t)\leq \gamma_{0}L(2H_{1}+1)\rho_{n}\left(
	\frac{1}{\nu}\frac{\omega(\ep)}{\ep}+
	\frac{1}{\sqrt{\nu}}\lk\omega(\ep)\right)E_{k,\ep}(t)
	\quad\quad\forall t\in[0,S].
	\label{est:ekep}
\end{equation}

Let us consider now the eigenvalues $\lk>\rho_{n}$, which are clearly
positive, and let us set $\ep_{k}:=\lk^{-1}$.  By (\ref{hp:phi-ndg}) we
have that 
$$\frac{\omega(\ep_{k})}{\ep_{k}}=\lk\omega(\ep_{k})=
\lk\omega\left(\frac{1}{\lk}\right)\leq\Lambda\varphi(\lk).$$

Using these estimates in (\ref{est:ekep}) we obtain that
$$E_{k,\ep_{k}}'(t)\leq \gamma_{0}L(2H_{1}+1)\rho_{n}\left(
\frac{1}{\nu}+\frac{1}{\sqrt{\nu}}\right)
\Lambda\varphi(\lk) E_{k,\ep_{k}}(t)=
\gamma_{2}\rho_{n}\varphi(\lk)E_{k,\ep_{k}}(t).$$

Integrating this differential inequality and using (\ref{defn:rhon})
we find that 
$$E_{k,\ep_{k}}(t)\leq E_{k,\ep_{k}}(0)
\exp\left(\gamma_{2}\rho_{n}\varphi(\lk)T\right)\leq E_{k,\ep_{k}}(0)
\exp\left(\rho_{n}^{2}\varphi(\lk)\right)$$
for every $t\in[0,S]$. Thanks to (\ref{est:cep-nu-mu}) we obtain that
\begin{eqnarray*}
	|u_{k}'(t)|^{2}+\lk^{2}|u_{k}(t)|^{2} & \leq & 
	\max\left\{1,\nu^{-1}\right\}E_{k,\ep_{k}}(t)\\
	 & \leq & \max\left\{1,\nu^{-1}\right\}
	 \left(|u_{1k}|^{2}+\lk^{2} c_{\ep_{k}}(0)|u_{0k}|^{2}\right)
	 \exp\left(\rho_{n}^{2}\varphi(\lk)\right)\\
	 & \leq & \max\left\{1,\nu^{-1}\right\}\cdot\max\{1,\mu\}
	 \left(|u_{1k}|^{2}+\lk^{2} |u_{0k}|^{2}\right)
	 \exp\left(\rho_{n}^{2}\varphi(\lk)\right)\\
	 & = & \gamma_{1}\left(|u_{1k}|^{2}+\lk^{2}|u_{0k}|^{2}\right)
	 \exp\left(\rho_{n}^{2}\varphi(\lk)\right),
\end{eqnarray*}
where $u_{0k}$ and $u_{1k}$ denote the Fourier components of $u_{0}$
and $u_{1}$, respectively.

By assumption (\ref{hp:data-ndg}) we have therefore that
$$\sum_{\lk>\rho_{n}}\lk\left(
|u_{k}'(t)|^{2}+\lk^{2}|u_{k}(t)|^{2}\right)\leq \gamma_{1}
\sum_{\lk>\rho_{n}}\lk \left(|u_{1k}|^{2}+\lk^{2}|u_{0k}|^{2}\right)
\exp\left(\rho_{n}^{2}\varphi(\lk)\right)\leq 2\gamma_{1}\rho_{n}$$
for every $t\in[0,S]$.  On the other hand, by (\ref{est:u'}) and
(\ref{est:au}) we have that
\begin{eqnarray*}
	\sum_{\lk\leq\rho_{n}}\lk\left(
	|u_{k}'(t)|^{2}+\lk^{2}|u_{k}(t)|^{2}\right) & \leq & 
	\rho_{n}\sum_{\lk\leq\rho_{n}}\left(
	|u_{k}'(t)|^{2}+\lk^{2}|u_{k}(t)|^{2}\right)\\
	 & \leq & \rho_{n}\left(|u'(t)|^{2}+|A^{1/2}u(t)|^{2}\right)  \\
	 & \leq & \rho_{n}\left(\mathcal{H}(0)+
	 \frac{\mathcal{H}(0)}{\nu}\right)
\end{eqnarray*}
for every $t\in[0,S]$. In particular for $t=S$ we have that
\begin{eqnarray*}
	\lefteqn{\hspace{-2em}\left|\langle A^{3/4}u(S),A^{1/4}u'(S) \rangle\right|  \leq  
	|A^{3/4}u(S)|^{2}+|A^{1/4}u'(S)|^{2}}  \\
	\noalign{\vspace{1ex}}
	\hspace{2em} & = & \sum_{\lk\leq\rho_{n}}\lk\left(
	|u_{k}'(S)|^{2}+\lk^{2}|u_{k}(S)|^{2}\right)+
	\sum_{\lk>\rho_{n}}\lk\left(
	|u_{k}'(S)|^{2}+\lk^{2}|u_{k}(S)|^{2}\right)  \\
	 & \leq & \rho_{n}\left(\mathcal{H}(0)+
	 \frac{\mathcal{H}(0)}{\nu}+2\gamma_{1}\right)  \\
	 & < & H_{1}\rho_{n}. 
\end{eqnarray*}

This contradicts (\ref{eq:S-nec}).

\subsection{The weakly hyperbolic case}

Let us introduce some constants.  Let $L$, $\Lambda$, $\gamma_{0}$ be
the constants appearing in (\ref{hp:m-ocont}), (\ref{hp:phi-dg}),
and in Lemma~\ref{lemma:conv}, and let
$$\gamma_{3}:=1+\frac{1}{\omega(1)},$$
$$\gamma_{4}:=\max\left\{\m{u(t)}:t\in[0,T/2]\right\}+
\max\left\{\omega(\sigma):0\leq \sigma\sqrt{\omega(\sigma)}\leq
1\right\},$$
$$\gamma_{5}:=\gamma_{3}(1+\gamma_{4})(\Lambda+1)$$
$$H_{2}:=\max\left\{\left|\langle A^{3/4}u_{0},A^{1/4}u_{1}
\rangle\right|+1,(|u_{0}|+1)\sqrt{\mathcal{H}(0)}+
\gamma_{5}+1)\right\},$$
$$\gamma_{6}:=1+\gamma_{0}L(2H_{2}+1).$$

Since $\rho_{n}\to +\infty$ we can choose $n\in\n$ such that
$\rho_{n}\geq 1$, and
\begin{equation}
	\rho_{n}^{1/2}\geq T\sqrt{\mathcal{H}(0)},
	\hspace{3em}
	\rho_{n}^{1/2}\geq 4\gamma_{6}\Lambda T,
	\hspace{3em}
	\rho_{n}\geq\frac{2}{T\sqrt{\omega(T/2)}}.
	\label{defn:rhon-w}
\end{equation}

Let us set 
$$S:=\sup\left\{\tau\leq T:\left|\langle A^{3/4}u(t),A^{1/4}u'(t)
\rangle\right|\leq H_{2}\rho_{n}^{5/2}\ \ \forall
t\in[0,\tau]\right\}.$$

We remark that $S>0$ because $\left|\langle
A^{3/4}u_{0},A^{1/4}u_{1} \rangle\right|<H_{2}\leq
H_{2}\rho_{n}^{5/2}$.  

If $S=T$ we can conclude as in the strictly hyperbolic case (using the
linear theory for the weakly hyperbolic case, for which we refer to
\cite{cjs}).  So let us assume that $S<T$.  By the maximality of $S$
we have that necessarily
\begin{equation}
	\left|\langle A^{3/4}u(S),A^{1/4}u'(S) \rangle\right|=
	H_{2}\rho_{n}^{5/2}.
	\label{eq:S-nec-w}
\end{equation}

Let us consider the function $c(t)$ defined according to
(\ref{defn:c}), let us extend it outside the interval $[0,S]$ as in
Lemma~\ref{lemma:conv}, and let us set
$$\cep(t):=\omega(\ep)+\int_{\re}^{}c(t+\ep s)\rho(s)\,ds
\quad\quad\forall t\in\re.$$

Arguing as in the strictly hyperbolic case we find that
$$\left|c(t)-c(s)\right| \leq L(2H_{2}+1)\rho_{n}^{5/2}\omega(|t-s|)$$ 
for every $t$ and $s$ in $[0,S]$.  Therefore from statement (4) of
Lemma~\ref{lemma:conv} we deduce that
\begin{equation}
	 |\cep(t)-c(t)|\leq\left(1+\gamma_{0}L(2H_{2}+1)
	 \rho_{n}^{5/2}\right)\omega(\ep)=
	 \gamma_{6}\rho_{n}^{5/2}\,\omega(\ep),  
	\label{est:cep-c-w}
\end{equation}
\begin{equation}
	|\cep'(t)|\leq\gamma_{0}L(2H_{2}+1)\rho_{n}^{5/2}
	\displaystyle{\frac{\omega(\ep)}{\ep}}\leq
	\gamma_{6}\rho_{n}^{5/2}\,
	\displaystyle{\frac{\omega(\ep)}{\ep}}.
	\label{est:cep'-w}
\end{equation}

Let us consider the Fourier components $u_{k}(t)$ of $u(t)$, and let
us define $E_{k,\ep}(t)$ as in (\ref{defn:ekep}).  Computing the time
derivative as in the strictly hyperbolic case, and using
(\ref{est:cep-c-w}), (\ref{est:cep'-w}), and the fact that
$\cep(t)\geq\omega(\ep)$ we find that
$$E_{k,\ep}'(t)\leq\gamma_{6}\rho_{n}^{5/2}\left(
\frac{1}{\ep}+\lk\sqrt{\omega(\ep)}\right)E_{k,\ep}(t)
\quad\quad\forall t\in[0,S].$$

Now we choose $\ep$ as a function of $k$.  The function
$h(\sigma)=\sigma\sqrt{\omega(\sigma)}$ is invertible.  Let us
consider the eigenvalues $\lk>\rho_{n}$, which are clearly positive,
and let us set $\ep_{k}:=h^{-1}(1/\lk)$.  By (\ref{hp:phi-dg}) we have
that
\begin{equation}
    \lk\sqrt{\omega(\ep_{k})}=\frac{1}{\ep_{k}}\leq\Lambda
    \varphi\left(\frac{1}{h(\ep_{k})}\right)=
    \Lambda\varphi(\lk),
    \label{est:Lambda}
\end{equation}
hence
$$E_{k,\ep_{k}}'(t)\leq 2\gamma_{6}\rho_{n}^{5/2}\Lambda
\varphi(\lk)E_{k,\ep_{k}}(t).$$

Integrating this differential inequality, and exploiting the second
condition in (\ref{defn:rhon-w}) we thus obtain that
$$E_{k,\ep_{k}}(t)\leq E_{k,\ep_{k}}(0)
\exp\left(2\rho_{n}^{5/2}\gamma_{6}\Lambda \varphi(\lk)T\right)\leq
E_{k,\ep_{k}}(0)
\exp\left(\frac{1}{2}\rho_{n}^{3}\varphi(\lk)\right)$$
for every $t\in[0,S]$.  In order to estimate $E_{k,\ep_{k}}(0)$ we
need an estimate on $c_{\ep_{k}}(0)$.  To this end we first observe
that $h(\ep_{k})=1/\lk<1$, hence
\begin{equation}
    \omega(\ep_{k})\leq\max\{\omega(\sigma):0\leq h(\sigma)\leq 1\}.
    \label{est:cep0-1}
\end{equation}

Moreover the last condition in (\ref{defn:rhon-w}) is equivalent to
$1/\rho_{n}\leq h(T/2)$.  Therefore from the monotonicity of $h$ it
follows that 
$$\ep_{k}=h^{-1}\left(\frac{1}{\lk}\right)\leq
h^{-1}\left(\frac{1}{\rho_{n}}\right)\leq
h^{-1}\left(h
\left(\frac{T}{2}\right)\right)=\frac{T}{2},$$
hence from statement (3) of Lemma~\ref{lemma:conv} we deduce that
\begin{equation}
    \int_{\re}c(\ep_{k}s)\rho(s)\,ds\leq
    \max\{c(t):0\leq t\leq\ep_{k}\}\leq
    \max\{c(t):0\leq t\leq T/2\}.
    \label{est:cep0-2}
\end{equation}

From (\ref{est:cep0-1}) and (\ref{est:cep0-2}) it follows that 
$c_{\ep_{k}}(0)\leq\gamma_{4}$, hence
$$E_{k,\ep_{k}}(0)\leq\max\left\{1,\cep(0)\right\}
\left(|u_{1k}|^{2}+\lk^{2}|u_{0k}|^{2}\right)\leq
(1+\gamma_{4})\left(|u_{1k}|^{2}+\lk^{2}|u_{0k}|^{2}\right).$$

Moreover from (\ref{th:omega-3}) and (\ref{est:Lambda}) it follows
that 
$$\max\left\{1,\frac{1}{\omega(\ep_{k})}\right\}\leq
1+\frac{1}{\omega(\ep_{k})}\leq
\gamma_{3}\left(1+\frac{1}{\ep_{k}}\right)\leq
\gamma_{3}(1+\Lambda\varphi(\lk)).$$

Since $(1+\Lambda x)\leq(\Lambda +1)e^{x/2}$ for every $\Lambda\geq 0$
and every $x\geq 0$, we have in particular that
$$\max\left\{1,\frac{1}{\omega(\ep_{k})}\right\}\leq
\gamma_{3}(1+\Lambda\varphi(\lk))\leq
\gamma_{3}(1+\Lambda)\exp\left(\frac{1}{2}\varphi(\lk)\right)\leq$$
$$\leq\gamma_{3}(1+\Lambda)\exp\left(\frac{1}{2}\rho_{n}^{3}
\varphi(\lk)\right).$$

From all these estimates it follows that
\begin{eqnarray*}
	|u_{k}'(t)|^{2}+\lk^{2}|u_{k}(t)|^{2} & \leq & 
	\max\left\{1,\frac{1}{\omega(\ep_{k})}\right\}E_{k,\ep_{k}}(t)\\
	 & \leq & \gamma_{3}(1+\Lambda)E_{k,\ep_{k}}(0)
	 \exp\left(\rho_{n}^{3}\varphi(\lk)\right)\\
	 & \leq & \gamma_{3}(1+\Lambda)(1+\gamma_{4})
	 \left(|u_{1k}|^{2}+\lk^{2}|u_{0k}|^{2}\right)
	 \exp\left(\rho_{n}^{3}\varphi(\lk)\right)  \\
	 & = & \gamma_{5} \left(|u_{1k}|^{2}+\lk^{2}|u_{0k}|^{2}\right)
	 \exp\left(\rho_{n}^{3}\varphi(\lk)\right).
\end{eqnarray*}

By assumption (\ref{hp:data-dg}) we have therefore that
$$\sum_{\lk>\rho_{n}}\lk\left(
|u_{k}'(t)|^{2}+\lk^{2}|u_{k}(t)|^{2}\right)\leq \gamma_{5}
\sum_{\lk>\rho_{n}}\lk \left(|u_{1k}|^{2}+
\lk^{2}|u_{0k}|^{2}\right)
\exp\left(\rho_{n}^{3}\varphi(\lk)\right)\leq 2\gamma_{5}\rho_{n}$$
for every $t\in[0,S]$, and in particular
\begin{eqnarray*}
    \left|\sum_{\lk>\rho_{n}}\lk^{2}u_{k}'(S)\cdot u_{k}(S)\right|
    & \leq & \sum_{\lk>\rho_{n}}\lk^{2}|u_{k}'(S)|\cdot|u_{k}(S)|
    \\
     & \leq & \frac{1}{2} \sum_{\lk>\rho_{n}}\left(
     \lk|u_{k}'(S)|^{2}+\lk^{3}|u_{k}(S)|^{2}\right) \\
     & \leq & \gamma_{5}\rho_{n}.
\end{eqnarray*}

On the other hand, by (\ref{est:u'}) and the
first condition in (\ref{defn:rhon-w}) we have that
$$|u(t)|\leq|u_{0}|+S\cdot\max\{|u'(t)|:t\in[0,S]\}
\leq|u_{0}|+T\cdot\sqrt{\mathcal{H}(0)}
\leq \left(|u_{0}|+1
\right)\rho_{n}^{1/2}$$
for every $t\in[0,S]$, hence
$$\left|\sum_{\lk\leq\rho_{n}}\lk^{2}u_{k}'(t) u_{k}(t)\right|
\leq\rho_{n}^{2}\left|\langle u(t),u'(t)\rangle\right|\leq
\rho_{n}^{2}|u(t)|\cdot|u'(t)| \leq \rho_{n}^{5/2}\left(|u_{0}|+1
\right)\sqrt{\mathcal{H}(0)}$$
for every $t\in[0,S]$. In particular for $t=S$ we have that
\begin{eqnarray*}
    \left|\langle A^{3/4}u(S),A^{1/4}u'(S) \rangle\right| & \leq &
    \left|\sum_{\lk\leq\rho_{n}}\lk^{2}u_{k}'(S)\cdot u_{k}(S)\right|+
    \left|\sum_{\lk>\rho_{n}}\lk^{2}u_{k}'(S)\cdot u_{k}(S)\right| \\
	\noalign{\vspace{1ex}}
     & \leq & \rho_{n}^{5/2}\left(|u_{0}|+1
	\right)\sqrt{\mathcal{H}(0)}+\gamma_{5}\rho_{n}  \\
	\noalign{\vspace{1ex}}
     & < & H_{2}\rho_{n}^{5/2}.
\end{eqnarray*}

This contradicts (\ref{eq:S-nec-w}).

\subsection{Proof of Proposition~\ref{prop:sum}}

Let us recursively define a sequence $\rho_{n}$ as follows. First of
all we set $\rho_{0}=0$. Let us assume that a term $\rho_{n}$ has been
defined. Assumption (\ref{hp:prop}) implies in particular that
$$ (u_{0},u_{1})\in
\G_{\varphi,r,3/4}(A)\times\G_{\varphi,r,1/4}(A)$$
with $r=\rho_{n}^{\beta}$, hence
$$\sum_{k=1}^{\infty}u_{0k}^{2}\lk^{3}\exp\left(
\rho_{n}^{\beta}\varphi(\lk)\right)<+\infty,
\hspace{3em}
\sum_{k=1}^{\infty}u_{1k}^{2}\lk\exp\left(
\rho_{n}^{\beta}\varphi(\lk)\right)<+\infty.$$

We can therefore choose $\rho_{n+1}$ big enough in such a way that
$\rho_{n+1}\geq\rho_{n}+1$, and
$$\sum_{\lk\geq\rho_{n+1}}u_{0k}^{2}\lk^{3}\exp\left(
\rho_{n}^{\beta}\varphi(\lk)\right)\leq\rho_{n},
\hspace{3em}
\sum_{\lk\geq\rho_{n+1}}^{\infty}u_{1k}^{2}\lk\exp\left(
\rho_{n}^{\beta}\varphi(\lk)\right)\leq\rho_{n}.$$

Let $\overline{u}_{0}$ and $\overline{u}_{1}$ be the elements of $H$
whose Fourier components are given by 
$$\overline{u}_{0k}:=\left\{
\begin{array}{ll}
	0 & \mbox{if }\rho_{2k}\leq\lk<\rho_{2k+1}, \\
	u_{0k} & \mbox{if }\rho_{2k+1}\leq\lk<\rho_{2k+2},
\end{array}
\right.
\hspace{1em}
\overline{u}_{1k}:=\left\{\begin{array}{ll}
	0 & \mbox{if }\rho_{2k}\leq\lk<\rho_{2k+1}, \\
	u_{1k} & \mbox{if }\rho_{2k+1}\leq\lk<\rho_{2k+2},
\end{array}
\right.$$
and let $\overline{\rho}_{n}:=\rho_{2n}$. We claim that
(\ref{th:prop-1}) holds true. Indeed for every $n\in\n$ we have that
\begin{eqnarray*}
	\sum_{\lk>\overline{\rho}_{n}}^{\infty}
	\overline{u}_{0k}^{2}\lk^{3}\exp\left(
	\overline{\rho}_{n}^{\beta}\varphi(\lk)\right) & = &
	\sum_{\lk>\rho_{2n}}^{\infty}\overline{u}_{0k}^{2}\lk^{3}
	\exp\left( \rho_{2n}^{\beta}\varphi(\lk)\right)\\
	 & = & \sum_{\lk\geq\rho_{2n+1}}^{\infty}\overline{u}_{0k}^{2}\lk^{3}
	 \exp\left( \rho_{2n}^{\beta}\varphi(\lk)\right) \\
	 & \leq & \sum_{\lk\geq\rho_{2n+1}}^{\infty}u_{0k}^{2}\lk^{3}
	 \exp\left( \rho_{2n}^{\beta}\varphi(\lk)\right) \\
	\noalign{\vspace{1ex}}
	 & \leq & \rho_{2n}=\overline{\rho}_{n},
\end{eqnarray*}
and similarly for $\overline{u}_{1}$. Note that in the second equality
we exploited the spectral gap of $\overline{u}_{0}$, whose 
components are equal to zero in the range $(\rho_{2n},\rho_{2n+1})$.

In the same way we can show that
$\widehat{u}_{0}:=u_{0}-\overline{u}_{0}$ and
$\widehat{u}_{1}:=u_{1}-\overline{u}_{1}$ satisfy (\ref{th:prop-2})
with $\widehat{\rho}_{n}:=\rho_{2n+1}$.
\qed

\label{NumeroPagine}

\end{document}